\documentclass[12pt]{amsart}
\usepackage{amssymb,amsfonts,amsmath,amsopn,amstext,amscd,color,latexsym,mathabx,mathrsfs,skak,stackrel,xy,url,verbatim}

\allowdisplaybreaks

\input xy
\xyoption{all}

\theoremstyle{remark}

\newtheorem{para}{\bf}[subsection]

\newtheorem{rem}[para]{\bf Remark}
\theoremstyle{definition}

\theoremstyle{plain}

\newtheorem{thm}[para]{Theorem}

\newtheorem{lemma}[para]{Lemma}

\newtheorem{prop}[para]{Proposition}

\newenvironment{numequation}{\addtocounter{para}{1}
\begin{equation}}{\end{equation}}

\newcommand{\bbG}{{\mathbb G}}

\newcommand{\bbP}{{\mathbb P}}

\newcommand{\bbX}{{\mathbb X}}

\newcommand{\bbZ}{{\mathbb Z}}

\newcommand{\frX}{{\mathfrak X}}

\newcommand{\cD}{{\mathcal D}}

\newcommand{\cF}{{\mathcal F}}

\newcommand{\cI}{{\mathcal I}}

\newcommand{\cO}{{\mathcal O}}

\newcommand{\cQ}{{\mathcal Q}}

\newcommand{\cT}{{\mathcal T}}

\newcommand{\sD}{{\mathscr D}}

\newcommand{\Q}{{\mathbb Q}}

\newcommand{\Z}{{\mathbb Z}}

\newcommand{\Qp}{{\mathbb Q_p}}

\newcommand{\Zp}{{\mathbb Z_p}}
\newcommand{\Fp}{{\mathbb F_p}}

\newcommand{\Pf}{{\it Proof. }}

\newcommand{\Spec}{{\rm Spec}}
\newcommand{\Spf}{{\rm Spf}}

\newcommand{\eqdef}{\;\stackrel{\text{\tiny def}}{=}\;}
\newcommand{\lra}{\longrightarrow}
\newcommand{\hra}{\hookrightarrow}

\newcommand{\midc}{{\; | \;}}
\newcommand{\ord}{{\rm ord}}
\newcommand{\ra}{\rightarrow}

\newcommand{\sub}{\subset}

\newcommand{\Sym}{{\rm Sym}}

\begin{document}

\title{Arithmetic differential operators on a semistable model of $\bbP^1$}
\author{Deepam Patel}
\address{Department of Mathematics, Purdue University,
150 N. University Street, West Lafayette, IN 47907, U.S.A.}
\email{deeppatel1981@gmail.com}
\author{Tobias Schmidt}
\address{Institut f\"ur Mathematik, Humboldt-Universit\"at zu Berlin,
Rudower Chaussee 25, 12489 Berlin, Germany}
\email{Tobias.Schmidt@mathematik.hu-berlin.de}
\author{Matthias Strauch}
\address{Indiana University, Department of Mathematics, Rawles Hall, Bloomington, IN 47405, U.S.A.}
\email{mstrauch@indiana.edu}

\thanks{D.P. would like to acknowledge support from the ANR program $p$-adic Hodge Theory and beyond (Th{\'e}HopaD) ANR-11-BS01-005. M.S. would like to acknowledge the support of the National Science Foundation (award DMS-1202303). T.S. would like to acknowledge support of the Heisenberg programme of Deutsche Forschungsgemeinschaft.}

\begin{abstract} In this paper we study sheaves of logarithmic arithmetic differential operators on a particular semistable model of the projective line.  
The main result here is that the first cohomology group of these sheaves is non-torsion. We also consider a refinement of the order filtration on the sheaf of level zero (before taking the $p$-adic completion). The associated graded sheaf, which we explicitly determine, explains to some extent the occurrence of the cohomology classes in degree one.
\end{abstract}

\maketitle

\tableofcontents

\section{Introduction}

In this paper we study sheaves of logarithmic arithmetic differential operators on a particular semistable model $\bbX_1$ of the projective line $\bbX = \bbX_0 = \bbP^1_\Zp$ over $\Zp$. This model is obtained by blowing up the reduced closed subscheme given by the set of $\Fp$-valued points of $\bbX$. We denote the corresponding formal schemes, the completions along the special fiber, by $\frX$ and $\frX_1$, respectively. The sheaf of logarithmic differential operators of level $m$, as defined in \cite[sec. 5]{PSS2}, will be denoted by $\cD^{(m)}_{\frX_1}$, and its $p$-adic completion by $\sD^{(m)}_{\frX_1}$. The formal scheme $\frX_1$ is the first member of a family of formal semistable models $\frX_n$ which we studied in \cite{PSS2}. In that paper, 
we obtained some results about the global sections of the sheaf of logarithmic arithmetic differential operators $\cD^{(m)}_{\frX_n}$. One question that had not been treated there was the relation between $H^0(\frX_n,\cD^{(m)}_{\frX_n})$ and $H^0(\frX_n,\sD^{(m)}_{\frX_n})$. More precisely, one may ask if the natural inclusion

$$\widehat{H}^0(\frX_n,\cD^{(m)}_{\frX_n}) \lra H^0(\frX_n,\sD^{(m)}_{\frX_n})$$

\vskip8pt

is an isomorphism. On the left hand side $\widehat{H}^0(\frX_n,\cD^{(m)}_{\frX_n})$ denotes the $p$-adic completion of $H^0(\frX_n,\cD^{(m)}_{\frX_n})$. It is straightforward to see that there is a canonical exact sequence

$$0 \ra \widehat{H}^0(\frX_n,\cD^{(m)}_{\frX_n}) \ra H^0(\frX_n, \sD^{(m)}_{\frX_n}) \ra T_p\left(H^1(\frX_n, \cD^{(m)}_{\frX_n})\right) \ra 0 \;,$$

\vskip8pt

where the group on the right is the $p$-adic Tate module

$$T_p\left(H^1(\frX_n, \cD^{(m)}_{\frX_n})\right) = \varprojlim_k H^1(\frX_n, \cD^{(m)}_{\frX_n})[p^k] \;,$$

\vskip8pt

of $H^1(\frX_n, \cD^{(m)}_{\frX_n})$. In this paper we only consider the case when $n=1$, and the main results are summarized in the following theorem.

\vskip8pt

{\bf Theorem.} {\it (i) $T_p\Big(H^1(\frX_1, \cD^{(m)}_{\frX_1})\Big) = 0$, and the map 

$$\widehat{H}^0(\frX_1,\cD^{(m)}_{\frX_1}) \lra H^0(\frX_1,\sD^{(m)}_{\frX_1})$$

\vskip8pt

is therefore an isomorphism.

\vskip5pt

(ii) There is a canonical surjective homomorphism

$$H^1(\frX_1, \sD^{(m)}_{\frX_1}) \lra \widehat{H}^1(\frX_1,\cD^{(m)}_{\frX_1}) \;,$$ 

\vskip8pt

and cohomology group on the right contains non-torsion elements. In particular, \linebreak $H^1(\frX_1,\sD^{(m)}_{\frX_1,\Q})$ does not vanish.

\vskip8pt

(iii) The cohomology group $H^1(\frX_1,\sD^\dagger_{\frX_1,\Q})$ does not vanish.}

\vskip8pt

The sheaf $\sD^\dagger_{\frX_1,\Q}$ in (iii) is the inductive limit of the sheaves $\sD^{(m)}_{\frX_1,\Q}$. 

\vskip8pt

The investigations here and in \cite{PSS2} were motivated by the question if the formal models $\frX_n$ mentioned above are $\sD^\dagger_{\frX_n,\Q}$-affine, and the non-vanishing of $H^1(\frX_1,\sD^\dagger_{\frX_1,\Q})$ gives therefore a negative answer when $n=1$. This has led us to consider in \cite{PSS4} a different family of sheaves $\widetilde{\sD}^{(m)}_{n,k,\Q}$ of $p$-adically complete differential operators on $\frX_n$, and as it is shown there, $\frX_n$ turns out to be $\widetilde{\sD}^{(m)}_{n,k,\Q}$-affine.

\vskip12pt

\section{Global sections and cohomology of $\sD^{(0)}$ on $\frX_1$}\label{level_zero}

Let $\bbX_1$ be the blow-up of the projective line $\bbX = \bbX_0 = \bbP^1_\Zp$ in the reduced closed subscheme given by the set of $\Fp$-valued points. We denote the corresponding formal schemes, the completions along the special fiber, by $\frX$ and $\frX_1$, respectively. For a more detailed discussion of this (formal) scheme we refer to \cite[sec. 4]{PSS2}.

\vskip8pt

\subsection{Cohomology groups and their completions}\label{completion_zero}

Let $\sD_{\frX_1} = \sD^{(0)}_{\frX_1}$ be the $p$-adic completion of the sheaf of logarithmic differential operators $\cD_{\bbX_1} = \cD_{\bbX_1}^{(0)}$, cf. \cite[sec. 5]{PSS2}. We write $\cD_{\frX_1}$ for the $\cO_{\frX_1}$-module generated by the restriction of $\cD_{\bbX_1}$ to $\frX_1$. 

\begin{lemma}\label{reduction} The canonical homomorphism

$$H^i(\frX_1, \sD_{\frX_1}) \lra  \varprojlim_k H^i(\frX_1, \cD_{\frX_1}/p^k\cD_{\frX_1})$$

\vskip8pt

is an isomorphism when $i=0$ and surjective if $i=1$. For $i>1$ source and target of this map vanish.
\end{lemma}

\Pf For an inverse system of sheaves $(\cF_k)_k$, the presheaf $U \mapsto \varprojlim_k \cF_k(U)$ is actually a sheaf. This gives the statement for $i=0$. For $i>1$ the source and target of the map vanish because $\frX_1$ is a noetherian topological space of dimension one. In order to treat the case $i=1$ we are going to use \cite[ch. 0, Prop. 13.3.1]{EGA_III}. The third condition of this proposition is fulfilled because the transition maps on the system of sheaves are obviously surjective. Let $U$ be an affine open subset of $\frX_1$. Denote by $\frX_{1,k}$ the reduction of $\frX_1$ modulo $p^k$, and let $U_k = U \times_{\frX_1} \frX_{1,k}$ be the open affine subset of $\frX_{1,k}$. Then we have for all $i>0$

$$H^i(U,\cD_{\frX_1}/p^k\cD_{\frX_1}) = H^i(U_k,\cD_{\frX_1}/p^k\cD_{\frX_1}) = 0 \;,$$

\vskip8pt

because $\cD_{\frX_1}/p^k\cD_{\frX_1}$ is a quasi-coherent sheaf on $\frX_{1,k}$. This shows that the second condition of loc.cit. is satisfied, and, for $i>0$, also the first condition. Consider the exact sequence of quasi-coherent sheaves on $\frX_{1,k+1}$

$$0 \lra p^k\cD_{\frX_1}/p^{k+1}\cD_{\frX_1} \lra \cD_{\frX_1}/p^{k+1}\cD_{\frX_1} \lra \cD_{\frX_1}/p^k\cD_{\frX_1} \lra 0 \;.$$

\vskip8pt

Because $p^k\cD_{\frX_1}/p^{k+1}\cD_{\frX_1}$ has vanishing first cohomology on $U_k$, this sequence stays exact after applying $H^0(U_k, - )$, and this shows that
the first condition of loc.cit. is fulfilled in the case $i=0$. Hence we can conclude that the map in question is surjective for $i=1$. \qed

\vskip8pt

Next we consider the tautological exact sequence of sheaves on $\frX_1$

$$0 \ra \cD_{\frX_1} \stackrel{p^k}{\lra} \cD_{\frX_1} \ra \cD_{\frX_1} /p^k\cD_{\frX_1} \ra 0 \;.$$

\vskip8pt

The long exact cohomology sequence to this sequence gives the exact sequence

$$H^i(\frX_1,\cD_{\frX_1})  \stackrel{p^k}{\lra}  H^i(\frX_1,\cD_{\frX_1}) \ra H^i(\frX_1,\cD_{\frX_1} /p^k\cD_{\frX_1}) \ra H^{i+1}(\frX_1,\cD_{\frX_1})  \stackrel{p^k}{\lra}  H^{i+1}(\frX_1,\cD_{\frX_1}) \;.$$

\vskip8pt

We thus get an exact sequence

\begin{numequation}\label{fund_exact_seq_k}
0 \ra  H^i(\frX_1,\cD_{\frX_1})\Big/p^kH^i(\frX_1,\cD_{\frX_1}) \ra H^i(\frX_1,\cD_{\frX_1} /p^k\cD_{\frX_1}) \ra H^{i+1}(\frX_1, \cD_{\frX_1})\left[p^k\right] \ra 0 \;,
\end{numequation}

where $H^{i+1}(\frX_1, \cD_{\frX_1})\left[p^k\right]$ denotes the subgroup of elements annihilated by multiplication by $p^k$. Put

$$\widehat{H}^i(\frX_1,\cD_{\frX_1}) = \varprojlim_k  \left(H^i(\frX_1,\cD_{\frX_1})\Big/p^kH^i(\frX_1,\cD_{\frX_1})\right) \;,$$

\vskip8pt

and

$$T_p\left(H^i(\frX_1, \cD_{\frX_1})\right) = \varprojlim_k H^i(\frX_1, \cD_{\frX_1})[p^k] \;,$$

\vskip8pt

where the transition map $H^i(\frX_1, \cD_{\frX_1})[p^k] \ra H^i(\frX_1, \cD_{\frX_1})[p^{k-1}]$ is the multiplication by $p$. We then have the 

\vskip8pt

\begin{prop}\label{fundexseq} (a) For all $i \ge 0$ there is a natural exact sequence

\begin{numequation}\label{fund_exact_seq}
0 \ra \widehat{H}^i(\frX_1,\cD_{\frX_1}) \ra \varprojlim_k H^i(\frX_1, \cD_{\frX_1}/p^k\cD_{\frX_1}) \ra T_p\left(H^{i+1}(\frX_1, \cD_{\frX_1})\right) \ra 0 \;.
\end{numequation}

(b) For $i=0$ the exact sequence in (a) becomes

\begin{numequation}\label{fund_exact_seq_H0}
0 \ra \widehat{H}^0(\frX_1,\cD_{\frX_1}) \ra H^0(\frX_1, \sD_{\frX_1}) \ra T_p\left(H^1(\frX_1, \cD_{\frX_1})\right) \ra 0 \;.
\end{numequation}

(c) The cohomology group $H^2(\frX_1, \cD_{\bbX_1})$ vanishes and the exact sequence in (a) gives therefore a canonical isomorphism

\begin{numequation}\label{fund_iso_H1}
\widehat{H}^1(\frX_1,\cD_{\frX_1}) \simeq \varprojlim_k H^1(\frX_1, \cD_{\frX_1}/p^k\cD_{\frX_1}) \;.
\end{numequation}

\end{prop}

\Pf (a) For varying $k$ the projective system 

$$H^i(\frX_1,\cD_{\frX_1})\Big/p^kH^i(\frX_1,\cD_{\frX_1})$$

\vskip8pt

has obviously surjective transition maps (hence satisfies the Mittag-Leffler condition). We can thus pass to the limit over $k$ and using \ref{reduction} we obtain the exact sequence \ref{fund_exact_seq}.

\vskip8pt

(b) We use (a) in the case $i=0$ and \ref{reduction}.

\vskip8pt

(c) $H^2(\frX_1, \cD_{\bbX_1})$ vanishes because $\frX_1$ is a noetherian space of dimension one. The stated isomorphism follows then directly from (a). \qed

\subsection{Vanishing of ${\rm R}^1{\rm pr}_*(\cD_{\frX_1})$}

We use the Leray spectral sequence for the blow-up morphism

$${\rm pr}: \frX_1 \lra \frX = \frX_0 \;.$$

\vskip8pt

Applied to the sheaf $\cD_{\frX_1}$ we get an exact sequence

\begin{numequation}\label{second_ex_seq}
0 \ra H^1(\frX, {\rm pr}_*(\cD_{\frX_1})) \ra H^1(\frX_1,\cD_{\frX_1}) \ra H^0(\frX,{\rm R}^1 {\rm pr}_*(\cD_{\frX_1})) \ra 0 \;.
\end{numequation}

Denote by $\cD_{\frX,d}$ and $\cD_{\frX_1,d}$ the sheaves of differential operators of degree less or equal to $d$.

\vskip8pt

\begin{lemma}\label{vanishing_R1} (a) For all $d \ge 0$ one has ${\rm R}^1 {\rm pr}_*(\cD_{\frX_1,d}) = 0$.

\vskip8pt

(b) ${\rm R}^1 {\rm pr}_*(\cD_{\frX_1}) = 0$.

\vskip8pt

(c) $H^1(\frX, {\rm pr}_*(\cD_{\frX_1})) = H^1(\frX_1,\cD_{\frX_1})$.

\vskip8pt

\end{lemma}

\Pf (a) {\it Reduction: passage to the graded sheaves.} We have

$$\cT_{\frX_1}^{\otimes d} = \cD_{\frX_1,d}/\cD_{\frX_1,d-1} \;,$$

\vskip8pt

and we consider the tautological exact sequence

\begin{numequation}\label{filtration}
0 \lra \cD_{\frX_1,d-1} \lra \cD_{\frX_1,d} \lra \cT_{\frX_1}^{\otimes d} \lra 0 \;.
\end{numequation}

For $d=0$ we have $\cD_{\frX_1,0} = \cT_{\frX_1}^{\otimes 0}  = \cO_{\frX_1}$. Therefore, if we show 

$${\rm R}^1{\rm pr}_*(\cT_{\frX_1}^{\otimes d}) = 0$$ 

\vskip8pt

for all $d \ge 0$, then we can argue by induction and get
${\rm R}^1{\rm pr}_*(\cD_{\frX_1,d}) = 0$ for all d. Using that taking higher direct images commutes with inductive limits we get

$${\rm R}^1{\rm pr}_*(\cD_{\frX_1}) = 0 \;.$$

\vskip8pt

{\it Working with local coordinates.} Over the complement of $pr^{-1}(\bbX(\Fp))$ the blow-up morphism is an isomorphism, and the stalk of the sheaf ${\rm R}^1 {\rm pr}_*(\cD_{\frX_1})$ vanishes thus outside $\frX(\Fp)$. Consider a point $P \in \frX(\Fp)$. We may assume that $P$ corresponds to the point given by the ideal $(x,p)$ of the ring

$$R = \Zp\langle x \rangle\left[\frac{1}{x^{p-1}-1}\right] \;.$$

\vskip8pt

Then $\Spf(R)$ is an open neighborhood of $P$ in $\frX$. Put

$$R' = \Zp\langle x,z \rangle\left[\frac{1}{x^{p-1}-1},\frac{1}{z^{p-1}-1}\right]\Big/(xz-p) \;,$$

\vskip8pt

and $R'' = \Zp\langle t \rangle$, and identify the open subsets $\Spf(R')\left[\frac{1}{z}\right] \sub \Spf(R')$ and $\Spf(R'')\left[\frac{1}{t}\right] \sub \Spf(R'')$ via the relation $zt=1$. Then 

$$pr^{-1}(\Spf(R)) = \Spf(R') \cup \Spf(R'')$$ 

\vskip8pt

is an open neighborhood of the fiber $pr^{-1}(P)$. To show that the stalk of ${\rm R}^1 {\rm pr}_*(\cT_{\frX_1}^{\otimes d})$ at $P$ vanishes it suffices to show that 

$$H^1(pr^{-1}(U), \cT_{\frX_1}^{\otimes d}) = 0$$

\vskip8pt

for all affine open subsets $U \sub \Spf(R) \sub \frX$ containing $P$. Identify $\Spf(R)$ with a closed subset of $\Spf(R')$. Then we have $pr^{-1}(U) = U \cup \Spf(R'')$. Hence it suffices to show that 

$$H^1(V \cup \Spf(R''), \cT_{\frX_1}^{\otimes d}) = 0$$

\vskip8pt

for all affine open subsets $V \sub \Spf(R') \sub \frX_1$ containing $P$ (which we also consider as a point of $\frX_1$).

\vskip8pt

{\it Using \v Cech cohomology.} For such a $V$ the open subset
$V \cup \Spf(R'')$ always contains $\Spf(R')\left[\frac{1}{z}\right] = \Spf(R'')\left[\frac{1}{t}\right]$ and we may thus assume $\Spf(R')\left[\frac{1}{z}\right] \sub V$. Then we have 

$$V \cap \Spf(R'') = \Spf(R')\left[\frac{1}{z}\right] = \Spf(R'')\left[\frac{1}{t}\right] \;.$$

\vskip8pt

Then $H^1(V \cup \Spf(R''), \cT_{\frX_1}^{\otimes d})$ is equal to
the cokernel of the map

$$H^0\left(V,\cT_{\frX_1}^{\otimes d}\right) \oplus H^0\left(\Spf(R''),\cT_{\frX_1}^{\otimes d}\right) \lra H^0\left(\Spf(R'')\left[\frac{1}{t}\right],\cT_{\frX_1}^{\otimes d}\right) \;.$$

\vskip8pt

which sends $(s_1,s_2)$ to the difference of these sections when restricted to $\Spf(R'')\left[\frac{1}{t}\right]$. Any element in 

$$H^0\left(\Spf(R'')\left[\frac{1}{t}\right],\cT_{\frX_1}^{\otimes d}\right)$$

\vskip8pt

has the form $\left(\sum_{i \in \bbZ} a_it^i\right)\partial_t^{\otimes d}$. The sum $\left(\sum_{i \ge 0} a_it^i\right)\partial_t^{\otimes d}$
clearly extends to a section over $\Spf(R'')$. Note that we have in $\cT_{\frX_1}^{\otimes d}$

$$\partial_t^{\otimes d} = (-z^2\partial_z)^{\otimes d} = (-1)^d z^{2d} \partial_z^{\otimes d}$$

\vskip8pt

and therefore 

$$\left(\sum_{i < 0} a_it^i\right)\partial_t^{\otimes d}  = (-1)^d \left(\sum_{i < 0} a_i z^{-i+d}\right) z^d\partial_z^{\otimes d} \;,$$

\vskip8pt

and this extends to a section over $V$. 

\vskip8pt

(b) This follows from (a) and the fact that the higher direct image functor commutes with inductive limits.

\vskip8pt

(c) This is an immediate consequence of (b) and \ref{second_ex_seq}. \qed

\vskip8pt

\subsection{The cohomology group $H^1(\frX,{\rm pr}_*(\cD_{\frX_1}))$}

Consider the exact sequence \ref{filtration} and the corresponding sequence of direct images on $\frX$

\begin{numequation}\label{fil_direct_im}
0 \lra {\rm pr}_*\Big(\cD_{\frX_1,d-1}\Big) \lra {\rm pr}_*\Big(\cD_{\frX_1,d}\Big) \lra {\rm pr}_*\Big(\cT_{\frX_1}^{\otimes d}\Big) \lra {\rm R}^1 {\rm pr}_*\Big(\cD_{\frX_1,d-1}\Big) = 0 \;,
\end{numequation}

where we have used \ref{vanishing_R1} (a). We have 

$$H^1(\frX,{\rm pr}_*(\cD_{\frX_1})) = \varinjlim_d H^1(\frX,{\rm pr}_*(\cD_{\frX_1,d})) \;.$$

\vskip8pt

Because $\cD_{\frX_1,d}$ is coherent and ${\rm pr}$ is projective, the sheaf ${\rm pr}_*(\cD_{\frX_1,d})$ is coherent and \linebreak $H^1(\frX,{\rm pr}_*(\cD_{\frX_1,d}))$ is thus a finitely generated $\Zp$-module. Since the corresponding cohomology group on the generic fiber (in the sense of rigid geometry) vanishes (by GAGA and \cite{BB81}), we see that $H^1(\frX,{\rm pr}_*(\cD_{\frX_1,d}))$ is annihilated by a finite power of $p$. (We will give below a more precise description of $H^1(\frX,{\rm pr}_*(\cD_{\frX_1,d}))$ which shows directly that it is annihilated by a finite power of $p$.) In the proof of theorem \ref{H1_direct_sum} we will need the following elementary 

\begin{lemma}\label{transformation} Let $x$, $y$ be the standard coordinates on $\bbP^1$ satisfying $xy=1$. Then we have $\partial_y = -x^2\partial_x$ and, more generally, for any $s \in \bbZ_{\ge 1}$

$$\partial_y^s = (-1)^s \sum_{t=1}^s a_{s,t}x^{s+t}\partial_x^t \;,$$

\vskip8pt

where for all $s \ge 1$ and $1 \le t \le s$

\begin{numequation}\label{coeffs_formula}
a_{s,t} = {s \choose t}\frac{(s-1)!}{(t-1)!} \;, \hskip10pt \mbox{in particular}\,, \hskip6pt a_{s,1} = s! \hskip6pt \mbox{ and } \hskip6pt a_{s,s} = 1 \;.
\end{numequation}

\end{lemma}

\Pf We prove this by induction on $s$. The formula holds obviously in the case $s=1$. Assuming the formula to be correct for a given $s$, we have

$$\begin{array}{lcl}
\partial_y^{s+1} & = & (-x^2\partial_x)(-1)^s \left(\sum_{t=1}^s a_{s,t}x^{s+t}\partial_x^t\right) \\
&&\\
& = & (-1)^{s+1}\sum_{t=1}^s \left(a_{s,t}x^2(x^{s+t}\partial_x + (s+t)x^{s+t-1})\partial_x^t\right)\\
&&\\
& = & (-1)^{s+1}\sum_{t=1}^s \left(a_{s,t}x^{s+t+2}\partial_x^{t+1} + a_{s,t}(s+t)x^{s+t+1}\partial_x^t\right)\\
&&\\
& = & (-1)^{s+1} \left(a_{s,1}(s+1)x^{s+2}\partial_x\right. \\
&&\\
&& \hfill + \left.[\sum_{t=2}^s \left(a_{s,t-1} + a_{s,t}(s+t)\right)x^{s+1+t}\partial_x^{t+1}] + a_{s,s}x^{2s+2}\partial_x^{s+1}\right)
\end{array}$$

\vskip8pt

Using \ref{coeffs_formula} we then get for $2 \le t \le s$:

$$\begin{array}{lcl}
a_{s,t-1} + a_{s,t}(s+t) & = & {s \choose t-1}\frac{(s-1)!}{(t-2)!} + {s \choose t}\frac{(s-1)!}{(t-1)!}(s+t) \\
&&\\
&=& \frac{s!}{(s-t+1)!(t-1)!}\frac{(s-1)!}{(t-2)!} + \frac{s!}{(s-t)!t!}\frac{(s-1)!(s+t)}{(t-1)!}\\
&&\\
&=& \frac{s!(s-1)!}{(t-1)!} \left[\frac{1}{(t-2)!(s-t+1)} + \frac{s+t}{(s-t)!t!}\right]\\
&&\\
&=& \frac{s!(s-1)!}{(t-1)!}\left[\frac{(t-1)t+(s-t+1)(s+t)}{t!(s-t+1)!}\right]\\
&&\\
&=& \frac{s!(s-1)!}{(t-1)!}\frac{s(s+1)}{t!(s+1-t)!} = {s+1 \choose t}\frac{s!}{(t-1)!} = a_{s+1,t}
\end{array}$$
 
\vskip8pt

And finally $a_{s,1}(s+1) = s!(s+1) = (s+1)! = a_{s+1,1}$. \qed

\vskip8pt 
\begin{thm}\label{H1_direct_sum} For all $d \ge 1$ the canonical map 

\begin{numequation}\label{H1_inj}
H^1(\frX,{\rm pr}_*(\cD_{\frX_1,d-1})) \ra H^1(\frX,{\rm pr}_*(\cD_{\frX_1,d}))
\end{numequation}

coming from the long exact cohomology sequence associated to
\ref{fil_direct_im} is injective and embeds $H^1(\frX,{\rm pr}_*(\cD_{\frX_1,d-1}))$ as a direct summand of $H^1(\frX,{\rm pr}_*(\cD_{\frX_1,d}))$. Therefore, there is a splitting:

\begin{numequation}\label{H1_split}
H^1(\frX,{\rm pr}_*(\cD_{\frX_1,d})) =  H^1(\frX,{\rm pr}_*(\cD_{\frX_1,d-1})) \oplus H^1(\frX,{\rm pr}_*(\cT_{\frX_1}^{\otimes d})) \;.
\end{numequation}

\end{thm}

\Pf (i) We start with some preliminary considerations. The sheaf ${\rm pr}_*(\cD_{\frX_1,d})$ (resp. ${\rm pr}_*(\cT_{\frX_1}^{\otimes d})$) is naturally a subsheaf of $\cD_{\frX,d}$ (resp. $\cT_\frX^{\otimes d}$), cf. \cite[5.2]{PSS2}, and we denote by $Q_{\le d}$ (resp. $Q_d$) the quotient sheaf. Consider the commutative diagram:

$$\begin{array}{ccccccccc}
0 & \ra & {\rm pr}_*(\cD_{\frX_1,d-1}) & \ra & \cD_{\frX,d-1} & \ra & Q_{\le d-1} & \ra & 0\\
&& \downarrow && \downarrow && \downarrow &&\\
0 & \ra & {\rm pr}_*(\cD_{\frX_1,d}) & \ra & \cD_{\frX,d} & \ra & Q_{\le d} & \ra & 0\\
&& \downarrow && \downarrow && \downarrow &&\\
0 & \ra & {\rm pr}_*(\cT_{\frX_1}^{\otimes d}) & \ra & \cT_{\frX}^{\otimes d} & \ra & Q_d & \ra & 0
\end{array}$$

\vskip8pt

where the horizontal sequences are the tautological exact sequences.
The corresponding long exact sequences give rise to the commutative diagram

\begin{numequation}\label{comm_dia}
\begin{array}{ccccccc}
&&&&&&\\
H^0(\frX,{\rm pr}_*(\cD_{\frX_1,d-1})) & \hra & H^0(\frX,\cD_{\frX,d-1}) & \ra & H^0(\frX,Q_{\le d-1}) & \twoheadrightarrow  & H^1(\frX,{\rm pr}_*(\cD_{\frX_1,d-1}))\\
\downarrow && \downarrow && \downarrow && \downarrow \\
H^0(\frX,{\rm pr}_*(\cD_{\frX_1,d})) & \hra & H^0(\frX,\cD_{\frX,d}) & \ra & H^0(\frX,Q_{\le d}) & \twoheadrightarrow & H^1(\frX,{\rm pr}_*(\cD_{\frX_1,d})) \\
\downarrow && \downarrow && \downarrow && \downarrow \\
H^0(\frX,{\rm pr}_*(\cT_{\frX_1}^{\otimes d})) & \hra & H^0(\frX,\cT_{\frX}^{\otimes d}) & \ra & H^0(\frX,Q_d) & \twoheadrightarrow & H^1(\frX,{\rm pr}_*(\cT_{\frX_1}^{\otimes d})) 
\end{array}
\end{numequation}

The sheaves $Q_{\le d-1}$, $Q_{\le d}$ and $Q_d$ are skyscraper sheaves with support in $\frX(\Fp)$.  Let $x_a$ be a local coordinate at $a \in \frX(\Fp)$. Then, cf. \cite[5.2 (c)]{PSS2},

\begin{numequation}\label{Q_le_d}
Q_{\le d} = \bigoplus_{a \in \frX(\Fp)} \bigoplus_{k=1}^{d} \bigoplus_{i=0}^{k-1} \left(\bbZ/p^{k-i}\right) \cdot x_a^i \partial_{x_a}^k \;,
\end{numequation}

\begin{numequation}\label{Q_le_d_minus_1}
Q_{\le d-1} = \bigoplus_{a \in \frX(\Fp)} \bigoplus_{k=1}^{d-1} \bigoplus_{i=0}^{k-1} \left(\bbZ/p^{k-i}\right) \cdot x_a^i \partial_{x_a}^k \;,
\end{numequation}

and 

\begin{numequation}\label{Qd}
Q_d = \bigoplus_{a \in \frX(\Fp)} \bigoplus_{i=0}^{d-1} \left(\bbZ/p^{d-i}\right) \cdot x_a^i \partial_{x_a}^d \;. 
\end{numequation}

Hence there is a splitting

\begin{numequation}\label{Q_split}
Q_{\le d} = Q_{\le d-1} \oplus Q_d \;.
\end{numequation}

\vskip8pt

We introduce the following notation and terminology. For a global section $\delta \in H^0(\frX,\cD_{\frX,d})$ we denote its image in $Q_{\le d}$ by $Q_{\le d}(\delta)$. The component of this element in $Q_d$, according to the splitting \ref{Q_split}, will be denoted by $Q_d(\delta)$, and we denote the components in $\bigoplus_{i=0}^{d-1} \left(\bbZ/p^{d-i}\right) \cdot x_a^i \partial_{x_a}^d$ corresponding to $a \in \frX(\Fp)$ by $Q_{d,a}(\delta)$. We call $Q_d(\delta)$ (resp. $Q_{\le d}(\delta)$) the {\it local data in degree $d$} (resp. {\it in degree less or equal to $d$}) of $\delta$. Similarly we call $Q_{d,a}(\delta)$ the {\it local data in degree $d$ at $a$} of $\delta$.   

\vskip8pt 

(ii) Now we prove the injectivity of the map \ref{H1_inj}. 
The injectivity of this map is equivalent, by the long exact cohomology sequence attached to \ref{fil_direct_im}, to the surjectivity of the map

\begin{numequation}\label{H0_surj}
H^0(\frX,{\rm pr}_*(\cD_{\frX_1,d})) \ra H^0(\frX,{\rm pr}_*(\cT_{\frX_1}^{\otimes d}))
\end{numequation}

which appears on the right hand side of \ref{comm_dia}. We are going to prove that \ref{H0_surj} is surjective as follows: consider $\delta_1 \in H^0(\frX,{\rm pr}_*(\cT_{\frX_1}^{\otimes d}))$ and let $\delta \in H^0(\frX,\cT_{\frX}^{\otimes d})$ be its image. Then $Q_d(\delta) = 0$. The crucial step is to lift $\delta$ to an element $\widetilde{\delta} \in H^0(\frX,\cD_{\frX,d})$ in such a way that $Q_{\le d}(\widetilde{\delta}) = 0$. This implies that $\widetilde{\delta}$ does in fact come from an element (necessarily unique) $\widetilde{\delta}_1 \in H^0(\frX,{\rm pr}_*(\cD_{\frX,d}))$ which is a preimage of $\delta_1$ under the map \ref{H0_surj}.

\vskip8pt

We let $x = x_0$ and $y = x_\infty$. Then $\delta \in H^0(\frX,\cT_\frX^{\otimes d})$ can be written as

$$\delta = \sum_{s=0}^{d-1}A_sy^s\partial_y^{\otimes d} + \sum_{s'=0}^{d} B_{s'}x^{s'}\partial_x^{\otimes d} \in H^0(\frX,\cT_{\frX}^{\otimes d}) \;.$$

\vskip8pt

Let us consider in detail what it means that $Q_d(\delta) = 0$. For instance, if we write $\delta$ in terms of $\partial_y$, we have to use the transformation formula (in $H^0(\frX,\cT_{\frX}^{\otimes d})$): $x^{s'}\partial_x^{\otimes d} = \pm y^{2d-s'}\partial_y^{\otimes d}$, and

$$\delta  = \sum_{s=0}^{d-1}A_sy^s\partial_y^{\otimes d} + \sum_{s'=0}^{d} (\pm B_{s'})y^{2d-s'}\partial_y^{\otimes d} \;.$$

Since $s' \le d$ we have $2d-s' \ge d$, we see that the vanishing of the local data of $\delta$ in degree $d$ at $\infty$ imposes the condition that $p^{d-s} | A_s$ for $0 \le s \le d-1$. Similarly we find $p^{d-s'} | B_{s'}$ for $0 \le s' \le d$. 

\vskip8pt

We are looking for a preimage $\widetilde{\delta} \in H^0(\frX,\cD_{\frX,d})$ of $\delta$ 
whose image in $H^0(\frX,Q_{\le d})$ vanishes. We start by taking as a candidate the element $\widetilde{\delta}_d$ which is given by the same formula as $\delta$, but now the summands are considered to be global sections of $\cD_{\frX,d}$, i.e., 

$$\widetilde{\delta}_d = \sum_{s=0}^{d-1}A_sy^s\partial_y^d + \sum_{s'=0}^{d} B_{s'}x^{s'}\partial_x^d \in H^0(\frX,\cD_{\frX,d}) \;.$$

\vskip8pt

(We write $\partial_x^{\otimes d}$ when we consider it as a section of $\cT_\frX^{\otimes d}$, and we write $\partial_x^d$ when we consider it as a section of $\cD_{\frX, d}$.) By \ref{transformation} this is indeed a global section of $\cD_{\frX,d}$.

\vskip8pt

The problem that we are facing now is this: while the local data of $\widetilde{\delta}_d$ in degree $d$ vanish (by assumption), it will in general not be the case that the local data of $\widetilde{\delta}_d$ in degree $<d$ vanish as well. Our aim is to modify $\widetilde{\delta}_d$ by adding a global section of $\cD_{\frX,d-1}$ to it, such that the difference has vanishing local data in all degrees, hence comes from an element in ${\rm pr}_*\cD_{\frX_1, d}$.

\vskip8pt

In order to do so, we determine the local data of $\widetilde{\delta}_d$ at infinity in all degrees. Using \ref{transformation} we write

$$\begin{array}{rl}
\widetilde{\delta}_d = & \sum_{s=0}^{d-1}A_sy^s\partial_y^d + \sum_{s'=0}^{d} B_{s'}x^{s'}\partial_x^d \\
&\\
= & \sum_{s=0}^{d-1}A_sy^s\partial_y^d + \sum_{s'=0}^{d} B_{s'}(-1)^d \left(\sum_{e=1}^d a_{d,e}y^{d+e-s'}\partial_y^e\right)\\
&\\
= & \sum_{s=0}^{d-1}A_sy^s\partial_y^d + (-1)^d \sum_{e=1}^d a_{d,e}\left(\sum_{s'=0}^{d} B_{s'}y^{d+e-s'}\right) \partial_y^e \;.\\
&\\
\end{array}$$

\vskip8pt

Because $d+e-s' \ge e$ the term $y^{d+e-s'}\partial_y^e$ does not contribute to local data at infinity. So, in fact, $\widetilde{\delta}_d$ has vanishing local data at infinity in all degrees less or equal to $d$.

\vskip8pt

Now we analyze the local data at points $a \in \frX(\Fp) \setminus \{\infty\} = \Fp$. Let $\xi_a \in \Zp$ be a lift of $a$. We use \ref{transformation} again and write

\begin{alignat*}{1}
\widetilde{\delta}_d = & \sum_{s=0}^{d-1}A_sy^s\partial_y^d + \sum_{s'=0}^{d} B_{s'}x^{s'}\partial_x^d \\
&\\
= & \sum_{s'=0}^d B_{s'}x^{s'}\partial_x^d + \sum_{s=0}^{d-1} A_s(-1)^d \left(\sum_{e=1}^d a_{d,e}x^{d+e-s}\partial_x^e\right)\\
&\\
= & \sum_{s'=0}^d B_{s'}x^{s'}\partial_x^d + (-1)^d \sum_{e=1}^d a_{d,e}\left(\sum_{s=0}^{d-1} A_s x^{d+e-s}\right) \partial_x^e\\
&\\
= & \left(\sum_{s'=0}^d B_{s'}x^{s'}+(-1)^d\sum_{s=0}^{d-1} A_s x^{2d-s}\right)\partial_x^d + (-1)^d \sum_{e=1}^{d-1} a_{d,e}\left(\sum_{s=0}^{d-1} A_s x^{d+e-s}\right) \partial_x^e\\
&\\
= & \left(\sum_{s'=0}^d B_{s'}x^{s'}+(-1)^d\sum_{s=0}^{d-1} A_s x^{2d-s}\right)\partial_x^d \\
&\\
& \hfill + (-1)^d \sum_{e=1}^{d-1} a_{d,e}\left(\sum_{s=0}^{d-1} A_s (x_a+\xi_a)^{d+e-s}\right) \partial_{x_a}^e\\
&\\
= & \left(\sum_{s'=0}^d B_{s'}x^{s'}+(-1)^d\sum_{s=0}^{d-1} A_s x^{2d-s}\right)\partial_x^d \\
&\\
& \hfill + (-1)^d \sum_{e=1}^{d-1} a_{d,e}\left(\sum_{s=0}^{d-1} \sum_{k=0}^{d+e-s} {d+e-s \choose k} \xi_a^k A_s x_a^{d+e-s-k}\partial_{x_a}^e \right) \;. \\
&\\
\end{alignat*}

\vskip8pt

The term ${d+e-s \choose k} \xi_a^k A_s x_a^{d+e-s-k}\partial_{x_a}^e$ gives a non-zero contribution to the local data at $a$ in degree $e$ only if $d+e-s-k < e$, i.e., $d < s+k$, and in this case the contribution is modulo $p^{e-(d+e-s-k)} = p^{s+k-d}$. Since $s+k-d \le s+(d+e-s)-d = e$ and because $p^{d-s} | A_s$ we find that the contribution of ${d+e-s \choose k} \xi_a^k A_s x_a^{d+e-s-k}\partial_{x_a}^e$ vanishes if $d-s \ge e$. So we only need to pay attention to those terms for which $d-s<e$ or, equivalently, $d-e <s$. Going back to the fourth line of the display above, we write

\begin{alignat*}{1}
\widetilde{\delta}_d = & \left(\sum_{s'=0}^d B_{s'}x^{s'}+(-1)^d\sum_{s=0}^{d-1} A_s x^{2d-s}\right)\partial_x^d + (-1)^d \sum_{e=1}^{d-1} a_{d,e}\left(\sum_{s=0}^{d-1} A_s x^{d+e-s}\right) \partial_x^e\\
&\\
= & \left(\sum_{s'=0}^d B_{s'}x^{s'}+(-1)^d\sum_{s=0}^{d-1} A_s x^{2d-s}\right)\partial_x^d + (-1)^d \sum_{e=1}^{d-1} a_{d,e}\left(\sum_{s=0}^{d-e} A_s x^{d+e-s}\right) \partial_x^e \\
&\\
& \hfill + (-1)^d \sum_{e=1}^{d-1} a_{d,e}\left(\sum_{d-e < s < d} A_s x^{d+e-s}\right) \partial_x^e \;.\\
&\\
\end{alignat*}

\vskip8pt

As mentioned above, the terms in $\sum_{e=1}^{d-1} a_{d,e}\left(\sum_{s=0}^{d-e} A_s x^{d+e-s}\right) \partial_x^e$ do not contribute to the local data in degrees less than $d$. Note that $\sum_{1 < s < d}A_sy^{s-1}\partial_y^{d-1}$ is a global section of $\cD_{\frX,d-1}$. Now consider 

$$\widetilde{\delta}_{d,d-1} \eqdef \widetilde{\delta}_d + a_{d,d-1}\left(\sum_{1 < s < d}A_sy^{s-1}\right)\partial_y^{d-1} \;.$$

\vskip8pt

Because $d-1-(s-1) = d-s$ and because $p^{d-s} | A_s$ this differential operator has vanishing local data at infinity in degree $d-1$ (and in degree $d$). We write $\widetilde{\delta}_{d,d-1}$ in terms of powers of $\partial_x$ and find:

\begin{alignat*}{1}
\widetilde{\delta}_{d,d-1} = & \left(\sum_{s'=0}^d B_{s'}x^{s'}+(-1)^d\sum_{s=0}^{d-1} A_s x^{2d-s}\right)\partial_x^d + (-1)^d \sum_{e=1}^{d-1} a_{d,e}\left(\sum_{s=0}^{d-e} A_s x^{d+e-s}\right) \partial_x^e\\
&\\
& + (-1)^d \sum_{e=1}^{d-1} a_{d,e}\left(\sum_{d-e < s < d} A_s x^{d+e-s}\right) \partial_x^e \\
&\\
& \hfill + (-1)^{d-1} a_{d,d-1}\sum_{1<s<d} A_s \left(\sum_{e=1}^{d-1} a_{d-1,e} x^{d-1+e-(s-1)}\partial_x^e\right)\\
&\\
= & \left(\sum_{s'=0}^d B_{s'}x^{s'}+(-1)^d\sum_{s=0}^{d-1} A_s x^{2d-s}\right)\partial_x^d + (-1)^d \sum_{e=1}^{d-1} a_{d,e}\left(\sum_{s=0}^{d-e} A_s x^{d+e-s}\right) \partial_x^e\\
&\\
& + (-1)^d \sum_{e=1}^{d-2} a_{d,e}\left(\sum_{d-e < s < d} A_s x^{d+e-s}\right) \partial_x^e + (-1)^d  a_{d,d-1}\left(\sum_{1 < s < d} A_s x^{2d-1-s}\right) \partial_x^{d-1} \\
&\\
& \hfill + (-1)^{d-1} a_{d,d-1}\sum_{1<s<d} A_s \left(\sum_{e=1}^{d-1} a_{d-1,e} x^{d+e-s}\partial_x^e\right)\\
&\\
= & \left(\sum_{s'=0}^d B_{s'}x^{s'}+(-1)^d\sum_{s=0}^{d-1} A_s x^{2d-s}\right)\partial_x^d + (-1)^d \sum_{e=1}^{d-1} a_{d,e}\left(\sum_{s=0}^{d-e} A_s x^{d+e-s}\right) \partial_x^e\\
&\\
& + (-1)^d \sum_{e=1}^{d-2} a_{d,e}\left(\sum_{d-e < s < d} A_s x^{d+e-s}\right) \partial_x^e + (-1)^d  a_{d,d-1}\left(\sum_{1 < s < d} A_s x^{2d-1-s}\right) \partial_x^{d-1} \\
&\\
&  + (-1)^{d-1} a_{d,d-1}\sum_{1<s<d} A_s \left(\sum_{e=1}^{d-2} a_{d-1,e} x^{d+e-s}\partial_x^e\right)\\
&\\
& \hfill + (-1)^{d-1} a_{d,d-1}\sum_{1<s<d}A_s a_{d-1,d-1} x^{2d-1-s}\partial_x^{d-1}\\
&\\
= & \left(\sum_{s'=0}^d B_{s'}x^{s'}+(-1)^d\sum_{s=0}^{d-1} A_s x^{2d-s}\right)\partial_x^d + (-1)^d \sum_{e=1}^{d-1} a_{d,e}\left(\sum_{s=0}^{d-e} A_s x^{d+e-s}\right) \partial_x^e\\
&\\
& + (-1)^d \sum_{e=1}^{d-2} a_{d,e}\left(\sum_{d-e < s < d} A_s x^{d+e-s}\right) \partial_x^e \\
&\\
&  + (-1)^{d-1} a_{d,d-1}\sum_{e=1}^{d-2} a_{d-1,e} \left(\sum_{1<s<d} A_s x^{d+e-s}\right)\partial_x^e \;.\\
\end{alignat*}

\vskip8pt

We therefore see that $\widetilde{\delta}_{d,d-1}$ has vanishing local data in degrees $d$ and $d-1$. As above, in the last sum $\sum_{e=1}^{d-2} a_{d-1,e} \left(\sum_{1<s<d} A_s x^{d+e-s}\right)\partial_x^e$ all those terms with $d-s \ge e$ do not contribute local data, so we write $\widetilde{\delta}_{d,d-1}$ as the sum of

$$\left(\sum_{s'=0}^d B_{s'}x^{s'}+(-1)^d\sum_{s=0}^{d-1} A_s x^{2d-s}\right)\partial_x^d$$

\vskip8pt

and

$$(-1)^d \sum_{e=1}^{d-1} a_{d,e}\left(\sum_{s=0}^{d-e} A_s x^{d+e-s}\right) \partial_x^e + (-1)^{d-1} a_{d,d-1}\sum_{e=1}^{d-2} a_{d-1,e}\left(\sum_{s=2}^{d-e} A_s x^{d+e-s}\right) \partial_x^e$$

\vskip8pt

and

$$(-1)^d \sum_{e=1}^{d-2} (a_{d,e}-a_{d,d-1}a_{d-1,e})\left(\sum_{d-e < s < d} A_s x^{d+e-s}\right) \partial_x^e$$\

\vskip8pt

Now we define

$$\widetilde{\delta}_{d,d-1,d-2} = \widetilde{\delta}_{d,d-1} - (-1)^d (a_{d,d-2}-a_{d,d-1}a_{d-1,d-2})\left( \sum_{2 < s < d} A_sy^{s-2}\right)\partial_y^{d-2} \;.$$

\vskip8pt

Continuing in this manner shows that we eventually find $\widetilde{\delta} \eqdef \widetilde{\delta}_{d,\ldots,1} \in H^0(\frX,\cD_{\frX,d})$ which has vanishing local data in all degrees less or equal to $d$, and its projection to $H^0(\frX,\cT_\frX^{\otimes d})$ is equal to $\delta$. This finishes the proof of the injectivity of the map \ref{H1_inj}.

\vskip8pt

(iii) Now we prove the splitting \ref{H1_split}. We start by making the following general remark: if $H_1$ is a subgroup of a finite abelian $p$-group $H$, then $H_1$ is a direct summand of $H$ if (and only if) $pH \cap H_1 = pH_1$. 

\vskip8pt

Now let $c_{\le d} \in H^0(\frX,Q_{\le d})$ be any element and let $[c_{\le d}] \in H^1(\frX,{\rm pr}_*\cD_{\frX,d})$ be its image. Suppose $p[c_{\le d}] = [pc_{\le d}]$ lies in $H^1(\frX,{\rm pr}_*\cD_{\frX,d-1})$, and write $[pc_{\le d}] = [c_{\le d-1}]$ for some element $c_{\le d-1} \in H^0(\frX,Q_{\le d-1})$. Then there is $\delta_{\le d} \in H^0(\frX,\cD_{\frX,d})$ such that $Q_{\le d}(\delta_{\le d}) = pc_{\le d}-c_{\le d-1}$. This is implies that 

$$Q_d(\delta_{\le d}) \in H^0(\frX,Q_d) = \bigoplus_{a \in \frX(\Fp)} \bigoplus_{i=0}^{d-1} \left(\bbZ/p^{d-i}\right) \cdot x_a^i \partial_{x_a}^d$$

\vskip8pt

is such that all its local data in the various groups $\bbZ/p^{d-i}$ are divisible by $p$. Write $\delta_{\le d} = \delta_d + \delta_{\le d-1}$ with 

$$\delta_d  = \sum_{s=0}^{d-1}A_sy^s\partial_y^d + \sum_{s'=0}^{d} B_{s'}x^{s'}\partial_x^d \;,$$

\vskip8pt

and with $\delta_{\le d-1} \in H^0(\frX,\cD_{\frX,d-1})$. The local data in degree $d$ of $\delta_{\le d}$ (or, equivalently, $\delta_d$) at $a=\infty$ and $a=0$ can be read off immediately from this expression for $\delta_d$ and it follows that all coefficients $A_s$, $0 \le s \le d-1$, and $B_{s'}$, $0 \le s' \le d$, are divisible by $p$. So we can write $\delta_d = p\delta_d'$, and hence $\delta_{\le d} = p \delta_d' + \delta_{\le d-1}$. We then have $Q_{\le d}(\delta_{\le d}) = pQ_{\le d}(\delta_d')+ Q_{\le d-1}(\delta_{\le d-1})$. From $Q_{\le d}(\delta_{\le d}) = pc_{\le d}-c_{\le d-1}$ we thus get 

$$p(Q_{\le d}(\delta_d')-c_{\le d}) = -Q_{\le d-1}(\delta_{\le d-1}) - c_{\le d-1} \;.$$

\vskip8pt

Write $Q_{\le d}(\delta_d')-c_{\le d} = c_d+c'_{\le d-1}$ with $c'_d \in Q_d$ and $c'_{\le d-1} \in Q_{\le d-1}$ and we find:

$$pc'_{\le d-1} = -Q_{\le d-1}(\delta_{\le d-1}) - c_{\le d-1} \;,$$

\vskip8pt

and thus $[c_{\le d-1}] = p[c'_{\le d-1}]$. \qed

\vskip12pt

In order to estimate the exponent of $H^1(\frX,{\rm pr}_*(\cD_{\frX_1,d}))$ we need the following elementary lemma.

\vskip12pt

\begin{lemma}\label{torsion_gp_quotient} Let $A = \bbZ/p^{n_1} \oplus \cdots \oplus \bbZ/p^{n_r}$ be an abelian torsion group with $0 < n_1 \le n_2 \le \ldots \le n_r$. Let $a \in A$ be an arbitrary element. Then $A/\langle a \rangle$ surjects onto $ \bbZ/p^{n_1} \oplus \cdots \oplus \bbZ/p^{n_{r-1}}$. 
\end{lemma}

\Pf Write $a = (a_1, \ldots, a_r)$, and choose $i \in \{1,\ldots,r\}$ such that 

$$\ord(a_i) = \max\{\ord(a_j) \midc j=1, \ldots,r\} \;.$$

\vskip8pt

If now $b = (b_1, \ldots,b_r) \in \langle a \rangle$ is such that $b_i = 0$, then $b=0$. Therefore, the map 

$$\bbZ/p^{n_1} \oplus \cdots \oplus \bbZ/p^{n_{i-1}} \oplus \bbZ/p^{n_{i+1}} \oplus \bbZ/p^{n_r} \hra \bbZ/p^{n_1} \oplus \cdots \oplus \bbZ/p^{n_r} = A \twoheadrightarrow A/\langle a \rangle$$

\vskip8pt

is injective. Because finite-abelian groups are self-dual (non-canonically), we see that there is a surjection 

$$A/\langle a \rangle \twoheadrightarrow \bbZ/p^{n_1} \oplus \cdots \oplus \bbZ/p^{n_{i-1}} \oplus \bbZ/p^{n_{i+1}} \oplus \bbZ/p^{n_r} \;.$$

\vskip8pt

But the group on the right clearly surjects onto $\bbZ/p^{n_1} \oplus \cdots \oplus \bbZ/p^{n_{r-1}}$. \qed

\vskip12pt

\begin{prop}\label{exist_of_torsion} For any $d \ge 1$ the cohomology group $H^1(\frX,{\rm pr}_*(\cT_{\frX_1}^{\otimes d}))$ contains elements of order $p^e$ where $e = \lfloor \frac{p-1}{p+1}(d+1) \rfloor$. In particular, as $d$ tends to infinity, the exponents of $H^1(\frX,{\rm pr}_*(\cT_{\frX_1}^{\otimes d}))$ and of $H^1(\frX,{\rm pr}_*(\cD_{\frX_1,d}))$ tend to infinity.
\end{prop}

\Pf By \ref{H1_split} we have 

$$H^1(\frX,{\rm pr}_*(\cD_{\frX_1,d})) =  H^1(\frX,{\rm pr}_*(\cD_{\frX_1,d-1})) \oplus H^1(\frX,{\rm pr}_*(\cT_{\frX_1}^{\otimes d})) \;.$$

\vskip8pt

Furthermore, $H^1(\frX,{\rm pr}_*(\cT_{\frX_1}^{\otimes d}))$ is the quotient of 

$$Q_d = \bigoplus_{a \in \frX(\Fp)} \bigoplus_{i=0}^{d-1} \left(\bbZ/p^{d-i}\right) \cdot x_a^i \partial_{x_a}^d \simeq  \left(\bigoplus_{i=0}^{d-1} \left(\bbZ/p^{d-i}\right)\right)^{\oplus (p+1)}\;,$$

\vskip8pt

cf. \ref{Qd}, by the image of $H^0(\frX,\cT_{\frX}^{\otimes d})$ which is a free $\Zp$-module of rank $2d+1$. Write $2d+1 = k(p+1)+r$ with $0 \le r \le p$, so that $k = \frac{2d+1}{p+1}-\frac{r}{p+1}$. Then, by applying \ref{torsion_gp_quotient} repeatedly we see that $H^1(\frX,{\rm pr}_*(\cT_{\frX_1}^{\otimes d}))$ must be of exponent at least $p^e$ where 

$$e = d-k = \frac{p-1}{p+1}d+\frac{r-1}{p+1} =  \frac{p-1}{p+1}(d+1)-\frac{p-r}{p+1} = \left\lfloor \frac{p-1}{p+1}(d+1) \right\rfloor \;.$$

\qed

\vskip10pt

\begin{rem} With some more work it should also be possible to explicitly determine the structure of $H^1(\frX,{\rm pr}_*(\cD_{\frX_1,d}))$. The filtration techniques in sec. \ref{sec_noetherian} might be helpful in doing so. The appendix \ref{H1_examples} contains some information on the structure of this cohomology group for small $p$ and $d$. 
\end{rem}

\vskip12pt

\subsection{$\widehat{H}^0(\frX_1,\cD_{\frX_1}) = H^0(\frX_1, \sD_{\frX_1})$ and $H^1(\frX_1, \sD_{\frX_1})$ is non-torsion}

\begin{thm}\label{main_I_zero} (a) $T_p H^1(\frX,{\rm pr}_*(\cD_{\frX_1})) = 0$.

\vskip8pt

(b) $T_p H^1(\frX_1,\cD_{\frX_1}) = 0$.

\vskip8pt

(c) $\widehat{H}^0(\frX_1,\cD_{\frX_1}) = H^0(\frX_1, \sD_{\frX_1})$.

\vskip8pt

(d) $\widehat{H}^1(\frX_1, \cD_{\frX_1})$ contains non-torsion elements.

\vskip8pt

(e) $H^1(\frX_1, \sD_{\frX_1})$ contains non-torsion elements.

\vskip8pt

\end{thm}

\Pf (a) We have $H^1(\frX,{\rm pr}_*(\cD_{\frX_1})) = \varinjlim_d H^1(\frX,{\rm pr}_*(\cD_{\frX_1,d}))$. Using \ref{H1_split} we see that

\begin{numequation}\label{inf_direct_sum}
H^1(\frX,{\rm pr}_*(\cD_{\frX_1})) \simeq \bigoplus_{d=1}^\infty H^1(\frX,{\rm pr}_*(\cT_{\frX_1}^{\otimes d})) \;.
\end{numequation}

Because each group $H^1(\frX,{\rm pr}_*(\cT_{\frX_1}^{\otimes d}))$ is a finite $p$-group, the $p$-adic Tate module 

$$T_p H^1(\frX,{\rm pr}_*(\cD_{\frX_1}))$$ 

\vskip8pt

must vanish.

\vskip8pt

(b) This follows from (a) and \ref{vanishing_R1}.

\vskip8pt

(c) This follows from (b) and \ref{fund_exact_seq_H0}.

\vskip8pt

(d) For $d \ge 1$ put $e_d = \lfloor \frac{p-1}{p+1}(d+1) \rfloor$. Let $c_d \in H^1(\frX,{\rm pr}_*(\cT_{\frX_1}^{\otimes d}))$ be an element of order $e_d$, cf. \ref{exist_of_torsion}. It follows from \ref{inf_direct_sum} that $H^1(\frX_1,\cD_{\frX_1}) = H^1(\frX,{\rm pr}_*(\cD_{\frX_1}))$ contains a subgroup isomorphic to $\bigoplus_{d \ge 1} \langle c_d \rangle$. Let $(n_d)_{d \ge 1}$ be an increasing sequence of non-negative integers $n_d \le e_d$ such that $\lim_{d \ra \infty} n_d = \infty$ and $\lim_{d \ra \infty} (e_d-n_d) = \infty$. Then $c = \sum_{d \ge 1} p^{n_d}c_d$ converges in the $p$-adic completion $\widehat{H}^1(\frX_1,\cD_{\frX_1})$ of $H^1(\frX_1,\cD_{\frX_1})$. Moreover, $c$ is clearly not a torsion element. 

\vskip8pt

(e) This follows from the fact that the map

$$H^1(\frX_1, \sD_{\frX_1}) \lra \varprojlim_k H^1(\frX_1, \cD_{\frX_1}/p^k\cD_{\frX_1}) = \widehat{H}^1(\frX_1,\cD_{\frX_1}) \;,$$

\vskip8pt

cf. \ref{reduction}, is surjective. (The equality sign on the right is \ref{fund_iso_H1}.) \qed

\vskip12pt

\section{Global sections and cohomology of $\sD^{(m)}$ on $\frX_1$}\label{level_m}

In this section we consider the sheaves of differential operators $\sD_{\frX_1}^{(m)}$ on $\frX_1$ of level $m \ge 0$. For their definition we refer to \cite{PSS2}. The discussion is along the same lines as in section \ref{level_zero}, with a few modifications which we are going to point out as we proceed.

\vskip8pt

\subsection{Comparing the cohomology of $\cD^{(m)}$ and $\sD^{(m)}$}\label{completion_m}

Let $\sD^{(m)}_{\frX_1}$ be the $p$-adic completion of the sheaf of logarithmic differential operators $\cD_{\bbX_1}^{(m)}$, considered as a sheaf on $\frX_1$. We write $\cD^{(m)}_{\frX_1}$ for the $\cO_{\frX_1}$-module generated by the restriction of $\cD^{(m)}_{\bbX_1}$ to $\frX_1$. The first lemma is exactly as \ref{reduction}.

\begin{lemma}\label{reduction_m} The canonical homomorphism

$$H^i(\frX_1, \sD^{(m)}_{\frX_1}) \lra \varprojlim_k H^i(\frX_1, \cD^{(m)}_{\frX_1}/p^k\cD^{(m)}_{\frX_1}) \;.$$ \qed

\vskip8pt

is an isomorphism for $i=0$ and is surjective for $i=1$. For $i>1$ source and target of this map vanish. 
\end{lemma}

And also the next result goes over without any changes.

\begin{prop}\label{fundexseq_m} (a) For all $i \ge 0$ there is a canonical exact sequence

\begin{numequation}\label{fund_exact_seq_m}
0 \ra \widehat{H}^i(\frX_1,\cD^{(m)}_{\frX_1}) \ra \varprojlim_k H^i(\frX_1, \cD^{(m)}_{\frX_1}/\cD^{(m)}_{\frX_1}) \ra T_p\left(H^{i+1}(\frX_1, \cD^{(m)}_{\frX_1})\right) \ra 0 \;.
\end{numequation}

(b) For $i=0$ the exact sequence in (a) is

\begin{numequation}\label{fund_exact_seq_H0_m}
0 \ra \widehat{H}^0(\frX_1,\cD^{(m)}_{\frX_1}) \ra H^0(\frX_1, \sD^{(m)}_{\frX_1}) \ra T_p\left(H^1(\frX_1, \cD^{(m)}_{\frX_1})\right) \ra 0 \;.
\end{numequation}

(c) The cohomology group $H^2\left(\frX_1, \cD^{(m)}_{\frX_1}\right)$ vanishes and the exact sequence in (a) gives therefore a canonical isomorphism

\begin{numequation}\label{fund_iso_H1_m}
\widehat{H}^1(\frX_1,\cD^{(m)}_{\frX_1}) \simeq \varprojlim_k H^1(\frX_1, \cD^{(m)}_{\frX_1}/p^k\cD^{(m)}_{\frX_1}) \;.
\end{numequation} \qed
\end{prop}

\vskip12pt

\subsection{Vanishing of ${\rm R}^1{\rm pr}_*(\cD^{(m)}_{\frX_1})$}

As above we use the Leray spectral sequence for the blow-up morphism

$${\rm pr}: \frX_1 \lra \frX = \frX_0 \;.$$

\vskip8pt

Applied to the sheaf $\cD^{(m)}_{\frX_1}$ we get an exact sequence

\begin{numequation}\label{second_ex_seq_m}
0 \ra H^1(\frX, {\rm pr}_*(\cD^{(m)}_{\frX_1})) \ra H^1(\frX_1,\cD^{(m)}_{\frX_1}) \ra H^0(\frX,{\rm R}^1 {\rm pr}_*(\cD^{(m)}_{\frX_1})) \ra 0 \;.
\end{numequation}

Denote by $\cD^{(m)}_{\frX,d}$ and $\cD^{(m)}_{\frX_1,d}$ the sheaves of differential operators of degree less or equal to $d$. Similar to \ref{filtration} we have an exact sequence

\begin{numequation}\label{filtration_m}
0 \lra \cD^{(m)}_{\frX_1,d-1} \lra \cD^{(m)}_{\frX_1,d} \lra (\cT_{\frX_1}^{\otimes d})^{(m)} \lra 0 \;.
\end{numequation}

For $m,d \ge 0$ we let $q^{(m)}_d$ be defined by $d = q^{(m)}_d p^m +r $ with $0 \le r < p^m$. In the proof of the lemma below we will use 

\begin{numequation}\label{divided_power_m}
(\cT_{\frX_1}^{\otimes d})^{(m)} =  \frac{q^{(m)}_d!}{d!} \cT_{\frX_1}^{\otimes d} \sub
\cT_{\frX_1}^{\otimes d} \otimes_\Zp \Qp \;,
\end{numequation}

cf. \cite[3.2]{PSS2}.

\begin{lemma}\label{vanishing_R1_m} (a) For all $d \ge 0$ one has ${\rm R}^1 {\rm pr}_*(\cD^{(m)}_{\frX_1,d}) = 0$.

\vskip8pt

(b) ${\rm R}^1 {\rm pr}_*(\cD^{(m)}_{\frX_1}) = 0$.

\vskip8pt

(c) $H^1(\frX, {\rm pr}_*(\cD^{(m)}_{\frX_1})) = H^1(\frX_1,\cD^{(m)}_{\frX_1})$. 

\vskip8pt
\end{lemma}

\Pf (a) This follows as in \ref{vanishing_R1} (a) using  \ref{divided_power_m} in the Cech cohomology argument.

\vskip8pt

(b) Follows from (a) by passing to the limit.

\vskip8pt

(c) Follows from (b) and \ref{second_ex_seq_m}.  \qed

\vskip12pt

\subsection{The cohomology group $H^1(\frX,{\rm pr}_*(\cD^{(m)}_{\frX_1}))$}

Consider the exact sequence \ref{filtration_m} and the corresponding sequence of direct images on $\frX$

\begin{numequation}\label{fil_direct_im_m}
0 \lra {\rm pr}_*\Big(\cD^{(m)}_{\frX_1,d-1}\Big) \lra {\rm pr}_*\Big(\cD^{(m)}_{\frX_1,d}\Big) \lra {\rm pr}_*\Big(\cT_{\frX_1}^{\otimes d}\Big) \lra {\rm R}^1 {\rm pr}_*\Big(\cD^{(m)}_{\frX_1,d-1}\Big) = 0 \;,
\end{numequation}

where we have used \ref{vanishing_R1_m} (a). We have 

$$H^1(\frX,{\rm pr}_*(\cD^{(m)}_{\frX_1})) = \varinjlim_d H^1(\frX,{\rm pr}_*(\cD^{(m)}_{\frX_1,d})) \;.$$

\vskip8pt

We put $\partial_x^{\langle d \rangle_{(m)}} =  \frac{q^{(m)}_d!}{d!} \partial_x^d$, and similarly for $\partial_{y}^d$ (and also for $\partial_{x_a}^d$). With this notation we deduce from \ref{transformation} the following

\begin{lemma}\label{transformation_m} Let $x$, $y$ be the standard coordinates on $\bbP^1$ satisfying $xy=1$. Then we have for any $s \in \bbZ_{\ge 1}$

$$\partial_y^{\langle s \rangle_{(m)}} = (-1)^s \sum_{t=1}^s a_{s,t}^{(m)} x^{s+t}\partial_x^{\langle t \rangle_{(m)}} \;,$$

\vskip8pt

where for all $s \ge 1$ and $1 \le t \le s$

\begin{numequation}\label{coeffs_formula_m}
a_{s,t}^{(m)} = {s \choose t}\frac{(s-1)!}{(t-1)!}\frac{q^{(m)}_s!}{s!}\left(\frac{q^{(m)}_t!}{t!}\right)^{-1} = {s-1 \choose t-1} \frac{q^{(m)}_s!}{q^{(m)}_t!} \;.
\end{numequation}

These numbers are always integers, and we have, in particular,

$$a_{s,1} = q^{(m)}_s! \hskip6pt \mbox{ and } \hskip6pt a_{s,s} = 1 \;.$$ \qed

\end{lemma}

\vskip8pt 

\begin{thm}\label{H1_direct_sum_m} For all $d \ge 1$ the canonical map 

\begin{numequation}\label{H1_inj_m}
H^1(\frX,{\rm pr}_*(\cD^{(m)}_{\frX_1,d-1})) \ra H^1(\frX,{\rm pr}_*(\cD^{(m)}_{\frX_1,d}))
\end{numequation}

coming from the long exact cohomology sequence associated to
\ref{fil_direct_im_m} is injective and embeds $H^1(\frX,{\rm pr}_*(\cD^{(m)}_{\frX_1,d-1}))$ as a direct summand of $H^1(\frX,{\rm pr}_*(\cD^{(m)}_{\frX_1,d}))$. Therefore, there is a splitting:

\begin{numequation}\label{H1_split_m}
H^1\left(\frX,{\rm pr}_*(\cD^{(m)}_{\frX_1,d})\right) =  H^1\left(\frX,{\rm pr}_*(\cD^{(m)}_{\frX_1,d-1})\right) \oplus H^1\left(\frX,{\rm pr}_*((\cT_{\frX_1}^{\otimes d})^{(m)})\right) \;.
\end{numequation} \qed

\end{thm}

\Pf The proof proceeds along the lines of \ref{H1_direct_sum} taking into account the following points.

\vskip8pt

(i) The skyscraper sheaf $Q_d^{(m)})$ (resp. $Q_{\le d}^{(m)}$)
is defined, similar as before, as the quotient of $(\cT_\frX^{\otimes d})^{(m)}$ (resp. $\cD_{\frX,d}^{(m)}$) by   ${\rm pr}_*((\cT_{\frX_1}^{\otimes d})^{(m)})$ (resp. ${\rm pr}_*(\cD_{\frX_1,d}^{(m)})$). By \ref{divided_power_m}, the sheaf $Q_{\le d}^{(m)}$ (resp. $Q_d^{(m)}$) is actually isomorphic to the sheaf $Q_{\le d}$ (resp. $Q_d$). 

\vskip8pt

(ii) The subtle part is the proof of the injectivity. As in the proof of \ref{H1_direct_sum} consider an element $\delta$ of $H^0\left(\frX,(\cT_\frX^{\otimes d})^{(m)}\right)$ whose image in the group $H^0(\frX,Q_d^{(m)})$ vanishes. Then we want to lift it to an element $\widetilde{\delta} \in H^0\left(\frX,\cD_{\frX,d}^{(m)}\right)$ such that the image of $\widetilde{\delta}$ in $H^0(\frX,Q_{\le d}^{(m)})$ vanishes. The discussion now proceeds along exactly the same lines as before. The difference is that one has to use the transformation formula in \ref{transformation_m}. This does not affect the arguments because the coefficients $a^{(m)}_{s,t}$ are integral. \qed

\vskip12pt

\begin{prop}\label{exist_of_torsion_m} For any $d \ge 1$ the cohomology group $H^1\left(\frX,{\rm pr}_*((\cT_{\frX_1}^{\otimes d})^{(m)})\right)$ contains elements of order $p^e$ where $e = \lfloor \frac{p-1}{p+1}(d+1) \rfloor$. In particular, as $d$ tends to infinity, the exponents of $H^1\left(\frX,{\rm pr}_*((\cT_{\frX_1}^{\otimes d})^{(m)})\right)$ and of $H^1\left(\frX,{\rm pr}_*(\cD_{\frX_1,d}^{(m)})\right)$ tend to infinity.
\end{prop}

\Pf The proof of \ref{exist_of_torsion} carries over to the case $m>0$. \qed

\vskip12pt

\begin{thm}\label{main_I_m} (a) $T_p H^1(\frX,{\rm pr}_*(\cD^{(m)}_{\frX_1})) = 0$.

\vskip8pt

(b) $T_p H^1(\frX_1,\cD^{(m)}_{\frX_1}) = 0$.

\vskip8pt

(c) $H^0(\frX_1,\cD^{(m)}_{\frX_1})^\wedge = H^0(\frX_1, \sD^{(m)}_{\frX_1})$.

\vskip8pt

(d) $H^1(\frX_1, \sD^{(m)}_{\frX_1})$ contains non-torsion elements.

\vskip8pt

\end{thm}

\Pf The proof of \ref{main_I_zero} carries over to the case when $m>0$. \qed

\vskip12pt

\subsection{$H^0(\frX_1,\cD_{\frX_1}^{(0)})$ is noetherian}\label{sec_noetherian}

In this section we continue our study of the ring of global sections of $\cD_{\frX_1}^{(0)}$. Similar results should  also hold for $H^0(\frX_1,\cD_{\frX_1}^{(m)})$. 

\vskip8pt

We first consider the graded ring with regard to the filtration given by the degree, or order, of the differential operators, i.e.,

$$F_d^{\deg}\left(H^0(\frX_1,\cD_{\frX_1}^{(0)})\right) = H^0(\frX_1,\cD_{\frX_1,d}^{(0)}) \;.$$

\vskip8pt

We denote by ${\rm gr}^{\deg}\left(H^0(\frX_1,\cD_{\frX_1}^{(0)})\right)$ the corresponding graded ring. 

\begin{prop}\label{graded_degree_fil} (i) The sequence

$$0 \ra H^0(\frX_1, \cD^{(0)}_{\frX_1,d-1}) \ra H^0(\frX_1,\cD^{(0)}_{\frX_1,d}) \ra H^0(\frX_1,\cT_{\frX_1}^{\otimes d}) \ra 0$$

\vskip8pt

is exact.

\vskip8pt

(ii) There is a canonical isomorphism

$${\rm gr}^{\deg}\left(H^0(\frX_1,\cD_{\frX_1}^{(0)})\right) = \bigoplus_{d \ge 0} H^0(\frX_1,\cT_{\frX_1}^{\otimes d}) = H^0(\frX_1,\Sym(\cT_{\frX_1})) \;.$$

\end{prop}

\Pf (i) This follows from the exact sequence \ref{fil_direct_im} and the injectivity of the map
in \ref{H1_inj}, cf. \ref{H1_direct_sum}.

\vskip8pt

(ii) This is an immediate consequence of (i). \qed

\vskip12pt

\begin{para}\label{filtration_on_Sym} We now consider the sheaf of algebras $\Sym(\cT_{\frX_1})$ on $\frX_1$, and similarly the sheaf of algebras $\Sym(\cT_{\bbX_1})$ on $\bbX_1$. We have obviously

$$H^0\left(\frX_1, \Sym(\cT_{\frX_1})\right) = \bigoplus_{d \ge 0} H^0(\frX_1,\cT_{\frX_1}^{\otimes d}) \;,$$

\vskip8pt

and  

$$H^0\left(\bbX_1, \Sym(\cT_{\bbX_1})\right) = \bigoplus_{d \ge 0} H^0(\bbX_1,\cT_{\bbX_1}^{\otimes d}) \;.$$

\vskip8pt 

We will be analyzing these sheaves of algebras by considering their direct images on $\frX = \frX_0$ and $\bbX = \bbX_0$, respectively. For the moment we will consider the algebraic case and then deduce the corresponding results for the sheaves on the formal schemes. 

\vskip8pt

To this end we recall the ideal sheaf $\cI_d = \cI_{1,d}$ on $\bbX_1 = \bbP^1_{\Zp}$ from \cite{PSS2}. It is locally defined by the ideal $(x_a,p)^d \subset \Zp[x_a]$, and one has 

$${\rm pr}_*(\cT_{\bbX_1}^{\otimes d}) = \cI_d \cT_{\bbX_0}^{\otimes d} \;.$$

\vskip8pt

We filter this latter sheaf as follows:

$$p^d\cT_{\bbX_0}^{\otimes d} \sub \ldots \sub
p^{d-i}\cI_i \cT_{\bbX_0}^{\otimes d} \sub \ldots \sub 
\cI_d \cT_{\bbX_0}^{\otimes d} \;.$$

\vskip8pt

Next we consider the following filtration $(\cF_i)_{i \ge 0}$, concentrated in non-negative degrees, on 

$${\rm pr}_*\left(\Sym(\cT_{\bbX_1})\right)  = \bigoplus_{d \ge 0} \cI_d \cT_\bbX^{\otimes d} \;.$$

We put

$$\cF_0 = \bigoplus_{d \ge 0} p^d\cT_\bbX^{\otimes d}$$
 
and for $i >0$ we set

$$\cF_i = \cF_{i-1} + \bigoplus_{d \ge i} p^{d-i}\cI_i\cT_\bbX^{\otimes d}  =  \bigoplus_{0 \le d < i} \cI_d\cT_\bbX^{\otimes d} \oplus  \bigoplus_{d \ge i} p^{d-i}\cI_i\cT_\bbX^{\otimes d} \;.$$

\vskip8pt
\end{para}

\begin{prop}\label{quot_sheaf} For $0 < i \le d$, the quotient of $p^{d-i}\cI_i\cT_\bbX^{\otimes d}$ by $p^{d-(i-1)}\cI_{i-1}\cT^{\otimes d}$ is canonically isomorphic to $\cO_{\bbP^1_\Fp}(2d-i(p+1))$, considered as a sheaf on $\bbP^1_\Zp = \bbX$.

\end{prop}

\Pf Because $\cI_i\cT_\bbX^{\otimes d}$ and $\cI_{i-1}\cT_\bbX^{\otimes d}$ are $p$-torsion free, we can write

\begin{numequation}\label{quot} p^{d-i}\cI_i\cT_\bbX^{\otimes d} / p^{d-(i-1)}\cI_{i-1}\cT_\bbX^{\otimes d} \; = \; p^{d-i}\left[\cI_i\cT_\bbX^{\otimes d} / p\cI_{i-1}\cT_\bbX^{\otimes d}\right] 
 \; \simeq \; \cI_i\cT_\bbX^{\otimes d} / p\cI_{i-1}\cT_\bbX^{\otimes d} \;.
\end{numequation}

Denote by $\cQ$ the quotient on the right of \ref{quot}. Because $\cI_i \sub \cI_{i-1}$ we get $p\cI_i \sub p\cI_{i-1}$, and this shows that $\cQ$ is $p$-torsion. Furthermore, locally the ideal sheaf $\cI_i$ (resp. $p\cI_{i-1}$) is defined by the ideal $(x_a,p)^i$ (resp. $p(x_a,p)^{i-1}$). Now consider $(x_a,p)^i/p(x_a,p)^{i-1}$ as an ideal in the quotient ring $\Zp[x_a]/(p(x_a,p)^{i-1})$. As an $\Fp[x_a]$-module it is naturally isomorphic to the ideal $(x_a^i) \sub \Fp[x_a]$:

$$(x_a,p)^i/p(x_a,p)^{i-1} \stackrel{\simeq}{\lra} (x_a^i)  \sub \Fp[x_a] \;, \;\; (f \mbox{ mod } p(x_a,p)^{i-1}) \mapsto (f \mbox{ mod } p) \;.$$

\vskip8pt

Thus we find that $\cQ$ is isomorphic to the product of $\cT_{\bbP^1_{\Fp}}^{\otimes d}$ and the ideal sheaf on $\bbP^1_{\Fp}$ whose divisor is $\sum_{a \in \bbP^1(\Fp)} i \cdot a$. This proves the claim. \qed

\vskip8pt

\begin{para}\label{graded_Sym_sheaf} It follows from \ref{quot_sheaf} that the graded sheaf of algebras on $\bbX$

$${\rm gr}^{\cF}\left({\rm pr}_*\left(\Sym(\cT_{\bbX_1})\right)\right)$$

\vskip8pt

is isomorphic to

\begin{numequation}\label{graded_Sym} \left[\bigoplus_{d \ge 0} p^d\cT_\bbX^{\otimes d}\right] \oplus \bigoplus_{i>0} \left[\bigoplus_{d \ge i} \mbox{ ''}p^{d-i}\mbox{''} \cdot \cO_{\bbP^1_\Fp}(2d-i(p+1)) \right] \;.
\end{numequation}

Here the factor ''$p^{d-i}$'' in front of $\cO_{\bbP^1_\Fp}(2d-i(p+1))$ is only a 'formal factor' which has its origin in \ref{quot}. We find it convenient to keep track of it.

\end{para}

\begin{para}\label{fil_H0_Sym} We now consider the induced filtration $(H^0(\cF_i))_i$ on 

$$H^0\left({\rm pr}_*(\Sym(\cT_{\bbX_1}))\right) = \bigoplus_{d \ge 0} H^0(\bbX,\cI_d \cT_\bbX^{\otimes d}) \;.$$

We thus have for $i \ge 0$ 

$$H^0\left(\bbX,\cF_i\right) =  \bigoplus_{0 \le d < i} H^0\left(\bbX,\cI_d\cT_\bbX^{\otimes d}\right) \oplus  \bigoplus_{d \ge i} H^0\left(\bbX,p^{d-i}\cI_i\cT_\bbX^{\otimes d}\right) \;.$$

\vskip8pt
\end{para}

\begin{prop}\label{H0_graded_Sym} (i) The canonical map

$${\rm gr}^{H^0(\cF)}\Big(H^0\Big(\bbX, {\rm pr}_*(\Sym(\cT_{\bbX_1}))\Big)\Big) \lra H^0\Big(\bbX,{\rm gr}^{\cF}\Big({\rm pr}_*\left(\Sym(\cT_{\bbX_1})\right)\Big)\Big)$$

\vskip8pt

is an isomorphism. 

\vskip8pt

(ii) The rings 

$${\rm gr}^{H^0(\cF)}\Big(H^0\Big(\bbX,{\rm pr}_*(\Sym(\cT_{\bbX_1}))\Big)\Big)$$

and 

$$\left[\bigoplus_{d \ge 0} H^0(\bbX,p^d\cT_\bbX^{\otimes d})\right] \oplus \bigoplus_{i>0} \left[\bigoplus_{2d \ge i(p+1)} \mbox{ ''}p^{d-i}\mbox{''} \cdot H^0(\bbX_\Fp, \cO_{\bbP^1_\Fp}(2d-i(p+1))) \right]$$

\vskip8pt

are canonically isomorphic as graded rings.

\vskip8pt
 
(iii) The ring

$${\rm gr}^{H^0(\cF)}\Big(H^0\Big(\bbX, {\rm pr}_*(\Sym(\cT_{\bbX_1}))\Big)\Big)$$

\vskip8pt

is noetherian.

\vskip8pt

(iv) The ring

$$H^0\Big(\bbX, {\rm pr}_*(\Sym(\cT_{\bbX_1}))\Big) = \bigoplus_{d \ge 0} H^0(\bbX,\cI_d\cT_\bbX^{\otimes d}) = \bigoplus_{d \ge 0} H^0(\bbX_1,\cT_{\bbX_1}^{\otimes d}) = H^0\Big(\bbX_1,\Sym(\cT_{\bbX_1})\Big)$$

\vskip8pt

is noetherian.

\vskip8pt

(v) $H^0(\bbX_1,\cD_{\bbX_1}^{(0)})$ and $H^0(\frX_1,\cD_{\frX_1}^{(0)})$ are noetherian rings.

\vskip8pt

(vi) $H^0(\frX_1,\sD_{\frX_1}^{(0)})$ and $H^0(\frX_1,\sD_{\frX_1,\Q}^{(0)})$  are noetherian rings.
\end{prop}

\Pf (i) For $0 < i \le d$ we put $\cQ = \cI_i\cT_\bbX^{\otimes d} / p\cI_{i-1}\cT^{\otimes d}$ and consider the tautological exact sequence of sheaves

\begin{numequation}\label{quot_seq} 0 \ra p\cI_{i-1}\cT^{\otimes d} \ra \cI_i\cT_\bbX^{\otimes d} \ra \cQ  \ra 0 \;.
\end{numequation}

The assertion in (i) is equivalent to saying that the corresponding sequence of global sections 

\begin{numequation}\label{ex_H0} 
0 \ra H^0(\bbX,p\cI_{i-1}\cT^{\otimes d}) \ra H^0(\bbX,\cI_i\cT_\bbX^{\otimes d}) \ra H^0(\bbX,\cQ) \ra 0
\end{numequation}

is exact too, cf. \ref{quot}.
By \ref{quot_sheaf} we have $\cQ \simeq \cO_{\bbP^1_\Fp}(2d-i(p+1))$. Therefore, if $2d -i(p+1) <0$ the sequence \ref{ex_H0} is trivially exact. Now suppose that $2d - i(p+1) \ge 0$. Under this assumption we will show that $H^1(\bbX,p\cI_{i-1}\cT^{\otimes d})$ vanishes. As the sheaves $p\cI_{i-1}\cT^{\otimes d}$ and $\cI_{i-1}\cT^{\otimes d}$ are isomorphic, we will work with the latter. Let $\cI_{i-1}'$ be the ideal sheaf defined locally by $(x_a^{i-1}) \sub \Zp[x_a]$. (We recall that the coordinate function $x_a$ vanishes at some lift of $a \in \bbP^1(\Fp)$ in $\bbP^1(\Zp)$). This is a subsheaf of $\cI_{i-1}$. Let $\cQ' = \cI_{i-1} \cT^{\otimes d}/\cI_{i-1}'\cT^{\otimes d}$ be the quotient. Away from divisor ${\rm div}(\cI'_1)$ this sheaf vanishes, and it is thus supported on ${\rm div}(\cI'_1)$, which is affine (it is isomorphic to a disjoint union of $p+1$ copies of $\Spec(\Zp)$). Because the cohomology of $\cQ'$ is the same as that on its support (cf. \cite[ch. III, 2.10
 ]{HartshorneA}), we conclude that $H^1(\bbX,\cQ') = 0$, as quasi-coherent sheaves on affine schemes have vanishing higher cohomology. Furthermore, we have $H^1(\bbX,\cI'_{i-1}\cT^{\otimes d}) = 0$, because $\cI'_{i-1}\cT^{\otimes d} \simeq \cO_\bbX(2d-(i-1)(p+1))$ and $2d-(i-1)(p+1) >0$. By the long exact cohomology sequence associated to 

$$0 \ra \cI_{i-1}'\cT^{\otimes d} \ra \cI_{i-1} \cT^{\otimes d} \ra \cQ' \ra 0$$

\vskip8pt

we can conclude that $H^1(\bbX,\cI_{i-1} \cT^{\otimes d}) = 0$. This shows that \ref{ex_H0} is also exact when $2d-i(p+1) \ge 0$, and this is what we had to show.

\vskip8pt

(ii) This follows immediately from (i), and \ref{graded_Sym}, and the observation that the sheaf $\cO_{\bbP^1_\Fp}(2d-i(p+1))$ has vanishing global sections when $2d-i(p+1)<0$.

\vskip8pt

(iii) We will assume for simplicity that $p>2$. (Simple variants of the following arguments should also cover the case $p=2$.) By (ii) the ring in question is isomorphic to

\begin{numequation}\label{explicit} \left[\bigoplus_{k \ge 0} p^k H^0(\bbX,\cO_\bbX(2k))\right] \oplus \bigoplus_{i>0} \left[\bigoplus_{k \ge 0} \mbox{ ''}p^{i\frac{p-1}{2}+k}\mbox{''} \cdot H^0(\bbX_\Fp, \cO_{\bbX_\Fp}(2k)) \right] \;.
\end{numequation}

(Here we have set $d = i\frac{p+1}{2}+k$.) Let $T$ be a generator of the direct summand in degree $(i,k) = (1,0)$. (This summand is a one-dimensional $\Fp$-space.) Furthermore, denote the summand

$$\bigoplus_{k \ge 0} p^k H^0(\bbX,\cO_\bbX(2k))$$

\vskip8pt

in degree $i=0$ by $R_0$, which is known to be a noetherian ring. Then it is a simple matter to check that the ring in \ref{explicit} is isomorphic to
$R_0[T]/(pT)$ and is thus a noetherian ring. 

\vskip8pt

(iv) This follows from the fact that the filtration $(H^0(\cF_i))_{i \ge 0}$ is in non-negative degrees, and because, by (iii), the corresponding graded ring is noetherian, cf. \cite[1.6.9]{MCR}.

\vskip8pt

(v) This follows from the fact that the degree filtration $(H^0(\cD_{\bbX_1,d}^{(0)}))_{d \ge 0}$ is in non-negative degrees, and that by \ref{graded_degree_fil} the corresponding graded ring is the one appearing in (iv), which is noetherian, cf. \cite[1.6.9]{MCR}. We also remark that the two rings in (v) are actually the same. 

\vskip8pt

(vi) This follows from (v) and \ref{main_I_zero} (c), cf. \cite[3.2.3, 3.4.0.1]{BerthelotDI}. \qed

\vskip8pt

\begin{rem} The filtration $(H^0(\cF_i))_i$ on 

$${\rm gr}^{\deg}\left(H^0(\frX_1,\cD_{\frX_1}^{(0)})\right) = H^0(\Sym(\cT_{\frX_1}))$$

\vskip8pt

gives rise to a filtration $\tilde{\cF}$ on $H^0(\frX_1,\cD_{\frX_1}^{(0)})$ which is a refinement of the filtration by degree, and which has the property that

$${\rm gr}^{\tilde{\cF}}\Big(H^0(\frX_1,\cD_{\frX_1}^{(0)})\Big) =  \left[\bigoplus_{k \ge 0} p^k H^0(\bbX,\cO_\bbX(2k))\right] \oplus \bigoplus_{i>0} \left[\bigoplus_{k \ge 0} \mbox{ ''}p^{i\frac{p-1}{2}+k}\mbox{''} \cdot H^0(\bbX_\Fp, \cO_{\bbX_\Fp}(2k)) \right] \;.$$

\vskip8pt

As we have shown in the proof of \ref{H0_graded_Sym} (iii), this ring is isomorphic to $R_0[T]/(pT)$ with a generator $T$ which is in degree $(i,k) = (1,0)$, cf. the proof of \ref{H0_graded_Sym} (iii) for the notation. The ''$d$-degree'' of $T$ is thus $d^* := \frac{p+1}{2}$. This gives rise to the following vague question. Is there an element $\delta \in H^0(\bbX_1,\cD_{\bbX_1,d^*})$ such that $H^0(\frX_1,\sD^\dagger_{\frX_1})$ is some kind of completion of $\cD^{an}(\bbG(1)^\circ)_{\theta_0}[\delta]$? (Cf. \cite{PSS2} for the relation of this distribution algebra to the global sections of the sheaf of differential operators.)
\end{rem}

\vskip12pt

\subsection{$H^1(\frX_1, \sD^\dagger_{\frX_1,\Q})$ does not vanish for $p>2$}

\begin{thm}\label{H1_dagger_not_vanishing} (i) The inductive limit 

$$\varinjlim_m \widehat{H}^1(\frX_1,\cD^{(m)}_{\frX_1}) \otimes_\Z \Q$$

\vskip8pt

does not vanish when $p>2$.

\vskip8pt

(ii) $H^1(\frX_1, \sD^\dagger_{\frX_1,\Q})$ does not vanish when $p > 2$.

\vskip8pt
\end{thm}

\Pf (i) Let us consider the transition map 

\begin{numequation}\label{transition}\widehat{H}^1(\frX_1, \cD^{(0)}_{\frX_1}) \lra \widehat{H}^1(\frX_1, \cD^{(m)}_{\frX_1}) \;.
\end{numequation}

Using \ref{inf_direct_sum} and its analogues in level $m$, together with \ref{vanishing_R1_m} and \ref{fund_iso_H1_m}, we rewrite 
\ref{transition} as 

\begin{numequation}\label{transition_explicit}
\left[\bigoplus_{d=1}^\infty H^1(\frX,{\rm pr}_*(\cT_{\frX_1}^{\otimes d}))\right]^\wedge \lra  \left[\bigoplus_{d=1}^\infty H^1\left(\frX,{\rm pr}_*((\cT_{\frX_1}^{\otimes d})^{(m)})\right)\right]^\wedge \;,
\end{numequation}

where $[-]^\wedge$ denotes the $p$-adic completion of $[-]$. Because of \ref{divided_power_m} we can formally write the right hand side of \ref{transition_explicit} as

$$\left[\bigoplus_{d=1}^\infty \frac{q^{(m)}_d!}{d!} H^1\left(\frX,{\rm pr}_*(\cT_{\frX_1}^{\otimes d})\right)\right]^\wedge \;,$$

\vskip8pt

and the map in \ref{transition_explicit} assumes the following explicit form

$$(c_d)_{d \ge 1} \mapsto \left(\frac{d!}{q^{(m)}_d!} \cdot c_d\right)_{d \ge 1} \;,$$

\vskip8pt

where $c_d \in H^1\left(\frX,{\rm pr}_*(\cT_{\frX_1}^{\otimes d})\right)$. 
Now let $c_d$ be a cohomology class of order $p^{e_d}$ where $e_d = 
\lfloor \frac{p-1}{p+1}(d+1) \rfloor$, cf. \ref{exist_of_torsion_m}. 
Denote by $v_p$ the (logarithmic) normalized $p$-adic valuation. Then we have

$$v_p\left(\frac{d!}{q^{(m)}_d!}\right) \le \frac{d}{p-1} - \left(\frac{\lfloor \frac{d}{p^m}\rfloor}{p-1} - \log_p\left(\left\lfloor \frac{d}{p^m}\right\rfloor\right)\right) \le \frac{d}{p-1} - \frac{d}{(p-1)p^m} + \log_p(d)+1\;.$$

\vskip8pt 

Let $n_d$ be a non-negative integer, and denote by $\ord$ the order of an element. Then

$$\begin{array}{rcl} v_p\left(\ord\left(\frac{d!}{q^{(m)}_d!} \cdot p^{n_d} c_d\right)\right) & \ge & \frac{p-1}{p+1}(d+1) - 1 - n_d - (\frac{d}{p-1} - \frac{d}{(p-1)p^m}+\log_p(d)+1) \\
&&\\
& = & \left(\frac{p-1}{p+1} - \frac{1}{p-1} + \frac{1}{(p-1)p^m}\right)d - n_d -\log_p(d)-2\\
&&\\
& = & \left(\frac{p^2-3p}{p^2-1} + \frac{1}{(p-1)p^m}\right)d - n_d-\log_p(d)-2 \;.
\end{array}$$

\vskip8pt

If $p\ge 3$ and if we put, for instance, $n_d = \lfloor \sqrt{d} \rfloor$, $d \ge 1$, then we have, for any $m$,  

$$\lim_{d \ra \infty} \left[\left(\frac{p^2-3p}{p^2-1} + \frac{1}{(p-1)p^m}\right)d - n_d -\log_p(d)-2 \right] = \infty \;.$$

\vskip8pt

This means that the sequence of elements $(p^{n_d}c_d)$  defines an element $c$ of $\widehat{H}^1(\frX_1, \cD^{(0)}_{\frX_1})$, and this element $c$ has the property that its image in any group $\widehat{H}^1(\frX_1, \cD^{(m)}_{\frX_1})$ is non-torsion. The image of $c$ in 

$$\varinjlim_m \widehat{H}^1(\frX_1, \cD^{(m)}_{\frX_1})$$

\vskip8pt

will also not be a torsion element.

\vskip8pt

(ii) Because the maps 

$$H^1(\frX_1, \sD^{(m)}_{\frX_1}) \lra \varprojlim_k H^1(\frX_1, \cD^{(m)}_{\frX_1}/p^k\cD^{(m)}_{\frX_1}) = \widehat{H}^1(\frX_1, \cD^{(m)}_{\frX_1})$$ 

are surjective, cf. \ref{reduction_m}, the same is true after tensoring with $\Q$ and taking the limit for $m \ra \infty$. \qed

\section{Appendix: computer calculations of $H^1(\frX_1, \cD_{\frX_1,d})$}\label{H1_examples}

In the following tables the entry ''$m \times p^n$'' signifies 
a direct summand of $H^1(\frX_1,\cD_{\frX_1,d})$ isomorphic to $(\bbZ/p^n)^{\oplus m}$. These tables were calculated using MAGMA.

\vskip12pt

{\it $p=2$.} 

\vskip8pt

$\begin{array}{r|rrrr}
d & \\
1 & 0 &&&\\
2 & 1 \times 2 &&&\\
3 & 3 \times 2 &&&\\
4 & 6 \times 2 &&&\\
5 & 9 \times 2 & 1 \times 2^2 &&\\
6 & 12 \times 2 & 3 \times 2^2 &&\\
7 & 15 \times 2 & 6 \times 2^2 &&\\
8 & 18 \times 2 & 9 \times 2^2 & 1 \times 2^3 & \\
9 & 21 \times 2 & 12 \times 2^2 & 3 \times 2^3 & \\
10 & 24 \times 2 & 15 \times 2^2 & 6 \times 2^3 & \\
11 & 27 \times 2 & 18 \times 2^2 & 9 \times 2^3 & 1 \times 2^4\\
\end{array}$

\vskip12pt

{\it $p=3$.} 

\vskip8pt

$\begin{array}{r|rrrrr}
d & \\
1 & 1 \times 3 &&&&\\
2 & 4 \times 3 &&&&\\
3 & 8 \times 3 & 1 \times 3^2 &&&\\
4 & 12 \times 3 & 4 \times 3^2 &&&\\
5 & 16 \times 3 & 8 \times 3^2 & 1 \times 3^3&&\\
6 & 20 \times 3 & 12 \times 3^2 &4 \times 3^3&&\\
7 & 24 \times 3 & 16 \times 3^2 &8 \times 3^3& 1 \times 3^4&\\
8 & 28 \times 3 & 10 \times 3^2 & 12 \times 3^3 & 4 \times 3^4&\\
9 & 32 \times 3 & 24 \times 3^2 & 16 \times 3^3 & 8 \times 3^4& 1 \times 3^5\\
\end{array}$

\vskip12pt

{\it $p=5$.} 

\vskip8pt

$\begin{array}{r|rrrrrr}
d & \\
1 & 3 \times 5 &&&&&\\
2 & 9 \times 5 & 1 \times 5^2 &&&&\\
3 & 15 \times 5 & 6 \times 5^2 &&&&\\
4 & 21 \times 5 & 12 \times 5^2 & 3 \times 5^3 &&&\\
5 & 27 \times 5 & 18 \times 5^2 & 9 \times 5^3 & 1 \times 5^4&&\\
6 & 33 \times 5 & 24 \times 5^2 & 15 \times 5^3 & 6 \times 5^4&&\\
7 & 39 \times 5 & 30 \times 5^2 & 21 \times 5^3 & 12 \times 5^4&3 \times 5^5&\\
8 & 45 \times 5 & 36 \times 5^2 & 27 \times 5^3 & 18 \times 5^4&9 \times 5^5& 1 \times 5^6\\
\end{array}$

\vskip12pt

{\it $p=7$.} 

\vskip8pt

$\begin{array}{r|rrrrr}
d & \\
1 & 5 \times 7 &&&&\\
2 & 13 \times 7 & 3 \times 7^2 &&&\\
3 & 21 \times 7 & 11 \times 7^2 & 1 \times 7^3&&\\
4 & 29 \times 7 & 19 \times 7^2 & 8 \times 7^3 &&\\
5 & 37 \times 7 & 27 \times 7^2 & 16 \times 7^3 & 5 \times 7^4&\\
6 & 45 \times 7 & 35 \times 7^2 & 24 \times 7^3 & 13 \times 7^4& 3 \times 7^5\\
\end{array}$

\vskip12pt

{\it $p=11$.} 

\vskip8pt

$\begin{array}{r|rrrrr}
d & \\
1 & 9 \times 11 &&&&\\
2 & 21 \times 11 & 7 \times 11^2 &&&\\
3 & 33 \times 11 & 19 \times 11^2 & 5 \times 11^3&&\\
4 & 45 \times 11 & 31 \times 11^2 & 17 \times 11^3 & 3 \times 11^4&\\
5 & 57 \times 11 & 43 \times 11^2 & 29 \times 11^3 & 15 \times 11^4& 1 \times 11^5\\
\end{array}$

\bibliographystyle{plain}
\bibliography{mybib}

\end{document}